\documentclass{article}

\usepackage{amsmath, amssymb, amsbsy}
\usepackage{booktabs}
\usepackage{centernot}
\usepackage{caption, subcaption}
\usepackage{framed}
\usepackage{graphicx, color}
\usepackage{makecell}
\usepackage{multirow}
\usepackage[sort&compress, square]{natbib}
\usepackage{xfrac}

\DeclareGraphicsExtensions{.pdf}

\newcommand{\given}[1][]{\ensuremath{#1\mid}}
\DeclareMathOperator{\HR}{HR}

\newcommand{\Labbe}{L'Abb\'{e} }

\DeclareMathOperator{\OR}{OR}
\newcommand{\pcrude}{\ensuremath{p^*_\text{crude}}}
\newcommand{\pmarg}{\ensuremath{p_\text{marg}}}
\DeclareMathOperator{\Prst}{Pr_\std}
\DeclareMathOperator{\RD}{RD}
\DeclareMathOperator{\RR}{RR}
\newcommand{\std}{\ensuremath{\text{std}}}
\newcommand{\twobytwo}{2$\times$2}

\begin{document}

\title{Rothman diagrams: the geometry of causal inference in epidemiology}
\author{Eben Kenah}
\maketitle

\begin{abstract}
    Here, we explain and illustrate a geometric perspective on causal inference in cohort studies that can help epidemiologists understand the role of standardization in causal inference as well as the distinctions between confounding, effect modification, and noncollapsibility.
    For simplicity, we focus on a binary exposure $X$, a binary outcome $D$, and a binary confounder $C$ that is not causally affected by $X$.
    Rothman diagrams plot risk in the unexposed on the x-axis and risk in the exposed on the y-axis.
    The crude risks define one point in the unit square, and the stratum-specific risks define two other points in the unit square.
    These three points can be used to identify confounding and effect modification, and we show briefly how these concepts generalize to confounders with more than two levels.
    We propose a simplified but equivalent definition of collapsibility in terms of standardization, and we show that a measure of association is collapsible if and only if all of its contour lines are straight.
    We illustrate these ideas using data from a study conducted in Newcastle upon Tyne, United Kingdom, where the causal effect of smoking on 20-year mortality was confounded by age.
    We conclude that causal inference should be taught using geometry before using regression models.
\end{abstract}

\textbf{Key Messages}
\begin{itemize}
  \item Rothman diagrams, where the risk of disease in the unexposed is on the x-axis and the risk in the exposed is on the y-axis, provide a geometric perspective on causal inference from which the distinctions between confounding, effect modification, and noncollapsibility can be seen clearly.
  \item There is confounding when the crude point is outside the convex hull of the stratum-specific points.
  \item Effect modification of a measure of association $M$ occurs when the stratum-specific points are on different contour lines of $M$.
  \item A measure of association is collapsible if and only if all of its contour lines are straight.
\end{itemize}

\section{Rothman diagrams}
The \emph{unit square} is the set of all ordered pairs $(x, y) \in \mathbb{R}^2$ where both $x$ and $y$ are in the interval $[0, 1]$ (i.e., the interval including zero, one, and all real numbers in between).
To represent the point $(x, y)$, we can plot $x$ on a horizontal axis (the x-axis) and $y$ on a vertical axis (the y-axis).
\citet{rothman1975pictorial} introduced a geometric perspective on causal inference where $x$ represented the risk of disease in the unexposed or untreated and $y$ represented the risk of disease in the exposed or treated.
\citet{labbe1987meta} introduced a similar plot for meta-analysis, and these \emph{\Labbe\ plots} were used by~\citet{richardson2017modeling} to discuss the estimation of risk differences and risk ratios.

I first heard about \Labbe\ plots from Thomas Richardson (Department of Statistics, University of Washington) at the Summer Institute for Statistics and Modeling in Infectious Diseases at the University of Washington in 2009, and I have used them since 2014 to teach standardization, confounding, effect modification, and collapsibility.
Here, I will explain and illustrate this approach, expanding on~\citet{rothman1975pictorial} in the light of more recent advances in causal inference~\citep{hernanrobins2023}.
Because our goal is causal inference rather than meta-analysis and because of historical precedence, we will call our pictures \emph{Rothman diagrams}.

As an illustration, we will use a beautiful example of confounding given by~\citet{appleton1996ignoring}.
It occurred in a cohort of women who participated in a study of thyroid and heart disease in Newcastle upon Tyne, United Kingdom in 1972-1974~\citep{tunbridge1977spectrum}.
Their smoking status was measured in the original study, and their 20-year survival status was measured in a follow-up study~\citep{vanderpump1995incidence}.
Table~\ref{tab:crude} shows the crude \twobytwo{} table for smoking and 20-year mortality, and Table~\ref{tab:crudemeas} shows the corresponding crude measures of association estimated using binomial generalized linear models (GLMs).
Smokers had lower 20-year mortality than nonsmokers, and this difference was statistically significant.

\begin{table}
  \centering
  \begin{tabular}{l|rr|r}
    \toprule
                & Dead  & Alive   & Total \\
    \midrule
    Smoker      & 139   & 443     & 582 \\
    Nonsmoker   & 230   & 502     & 732 \\
    \midrule
    Total       & 369   & 945     & 1,314 \\
    \bottomrule
  \end{tabular}
  \caption{Crude \twobytwo{} table for smoking status in the original study and 20-year mortality~\citep{appleton1996ignoring}.}
  \label{tab:crude}
\end{table}

\begin{table}
  \centering
  \begin{tabular}{lrcr}
    \toprule
    Measure of              &           & Likelihood ratio (LR)     & \multicolumn{1}{c}{LR} \\
    association             & Estimate  & 95\% confidence interval  & p-value \\
    \midrule
    Odds ratio              & 0.685     & (0.535, 0.875)            & \multirow{4}{*}{0.0024} \\
    Risk ratio              & 0.760     & (0.633, 0.908)            & \\
    Risk difference         & -0.075    & (-0.123, -0.027)          & \\
    Hazard ratio            & 0.724     & (0.584, 0.892)            & \\
    \bottomrule
  \end{tabular}
  \caption{
    Crude measures of association for smoking and 20-year mortality with likelihood ratio confidence limits and p-value (which does not depend on the measure of association).
    Binomial GLMs used the logit link for the odds ratio, the log link for the risk ratio, the identity link for the risk difference, and the complementary log-log link for the hazard ratio.
  }
  \label{tab:crudemeas}
\end{table}

Table~\ref{tab:strat} shows \twobytwo\ tables stratified by age at the time of the original survey.
Older participants were less likely to smoke but more likely to die within 20 years than younger participants.
Among those aged 18-64 years, the prevalence of smoking was $\sfrac{533}{1{,}072} \approx 0.497$ and the risk of death was $\sfrac{162}{1{,}072} \approx 0.151$.
Among those aged $\geq 65$ years, the prevalence of smoking was $\sfrac{49}{242} \approx 0.202$ and the risk of death was $\sfrac{207}{242} \approx 0.855$.
The difference in survival between smokers and nonsmokers was due in part to different age distributions, so the causal effect of smoking on 20-year mortality was confounded by age.

\begin{table}
  \begin{subtable}{0.5\textwidth}
    \centering
    \begin{tabular}{l|rr|r}
      \toprule
      \multicolumn{4}{c}{Participants aged 18-64 years} \\
                  & Dead  & Alive   & Total \\
      \midrule
      Smoker      & 97    & 436     & 533 \\
      Nonsmoker   & 65    & 474     & 539 \\
      \midrule
      Total       & 162   & 910     & 1,072 \\
      \bottomrule
    \end{tabular}
  \end{subtable}
  \begin{subtable}{0.5\textwidth}
    \centering
    \begin{tabular}{l|rr|r}
      \toprule
      \multicolumn{4}{c}{Participants aged $\geq 65$ years} \\
                  & Dead  & Alive   & Total \\
      \midrule
      Smoker      & 42    & 7       & 49 \\
      Nonsmoker   & 165   & 28      & 193 \\
      \midrule
      Total       & 207   & 35      & 242 \\
      \bottomrule
    \end{tabular}
  \end{subtable}
  \caption{
    Age-stratified \twobytwo{} tables for smoking and 20-year mortality adapted from~\citet{appleton1996ignoring}.
    Ages are those at the time of participation in the original survey in 1972-1974.}
  \label{tab:strat}
\end{table}

\section{Risks and points}
Let $X$ be a binary exposure or treatment, $D$ be a binary disease outcome, and $C$ be a binary covariate that is not causally affected by $X$.
Let $D^x$ be the potential disease outcome when we intervene to set $X = x$~\citep{rubin1974estimating}.
The exposure, outcome, and covariate of individual $i$ are $X_i$, $D_i$, and $C_i$, respectively.
The potential outcome of $i$ under an intervention that sets $X_i = x$ is $D^x_i$.
If $X_i = x$, then $D_i = D^x_i$ by \emph{consistency} of potential outcomes~\citep{pearl2009causality}.

In our example, $X$ is smoking status at the time of the original survey, $D$ is death within 20 years, and $C$ is age.
The potential outcome $D^1$ is death within 20 years as a smoker, and $D^0$ is death within 20 years as a nonsmoker.
By consistency, $D = D^1$ for smokers and $D = D^0$ for nonsmokers.
We see only one of the two potential outcomes for each individual.

Probabilities are defined by proportions of the underlying population.
For example, $\Pr(X = x, C = c)$ is the proportion of the population with $X = x$ and $C = c$.
Conditional probabilities are defined by proportions of subpopulations.
For example, $\Pr(X = x \given C = c)$ is the proportion of the subpopulation with $C = c$ that has $X = x$.
Probabilities and conditional probabilities can estimated using a sample of $n$ individuals from the population, with uncertainty due to sampling variation summarized in terms of frequentist confidence intervals or Bayesian credible intervals.
In our example, the study participants can be seen as a sample of size $n = 1{,}314$ from the population of women in the United Kingdom between 1972 and 1974.
For simplicity, we will assume no selection bias and ignore sampling variation whenever possible.

\subsection{Counterfactual risks and causal points}
The counterfactual risk of disease when we intervene to set $X = 0$ in the entire population is $\Pr(D^0 = 1)$, and the counterfactual risk of disease when we set $X = 1$ is $\Pr(D^1 = 1)$.
The point
\begin{equation}
  \pmarg = \Big(\Pr(D^0 = 1), \Pr(D^1 = 1)\Big)
  \label{eq:mcpoint}
\end{equation}
on the Rothman diagram is the \emph{marginal causal point} for the population.
It is causal because it is based on counterfactual risks of disease under intervention, and it is marginal because it ignores the covariate $C$.

In the subpopulation with $C = c$, the counterfactual risk of disease when we set $X = x$ is:
\begin{equation}
  \Pr(D^x = 1 \given C = c) = \frac{\Pr(D^x = 1, C = c)}{\Pr(C = c)}
\end{equation}
The point
\begin{equation}
  p_c = \Big(\Pr\bigl(D^0 = 1 \given C = c\bigr), \Pr\bigl(D^1 = 1 \given C = c\bigr)\Big)
\end{equation}
on the Rothman diagram is the \emph{stratum-specific causal point} for the subpopulation with $C = c$.

In practice, we are almost always missing the data needed to plot causal points on a Rothman diagram.
To draw marginal or stratum-specific causal points, we need both potential outcomes $D^1_i$ and $D^0_i$ for each individual $i$ in our sample.
In our example, we would need to know whether the smokers would have survived for 20 years had they been nonsmokers and whether the nonsmokers would have survived for 20 years had they been smokers.

\subsection{Conditional risks and association points}
\begin{figure}
  \includegraphics[width = \textwidth]{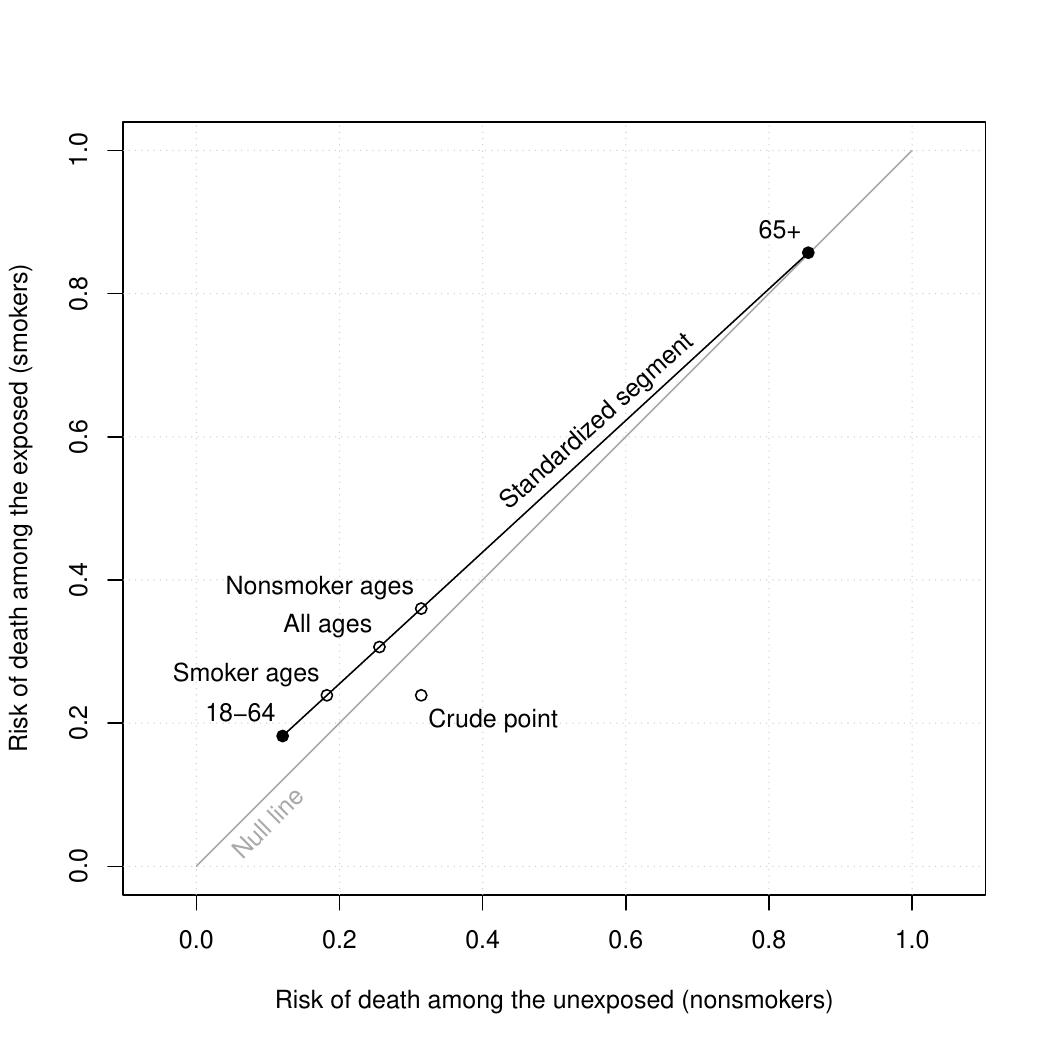}
  \caption{
    A Rothman diagram with the crude association point from Table~\ref{tab:crude} (open circle), the age-stratified association points from Table~\ref{tab:strat} (solid circles), the standardized segment (black line), and three standardized points (open circles).
    The null line is the gray diagonal line from the bottom left to the top right.
  }
  \label{fig:Rothman_std}
\end{figure}

Unlike causal points, we can plot \emph{association points} based on observed data.
The risk of disease among the unexposed is $\Pr(D = 1 \given X = 0)$, and the risk of disease among the exposed is $\Pr(D = 1 \given X = 1)$.
The point
\begin{equation}
  \pcrude = \Big(\Pr(D = 1 \given X = 0), \Pr(D = 1 \given X = 1)\Big)
\end{equation}
on the Rothman diagram is the \emph{crude association point} for the population, which is approximated by the crude association point for the study sample.
Figure~\ref{fig:Rothman_std} shows the crude point for the data from Table~\ref{tab:crude}.
Its x-coordinate is the risk of death among nonsmokers, which is $\sfrac{230}{702} \approx 0.314$.
Its y-coordinate is the risk of death among smokers, which is $\sfrac{139}{582} \approx 0.239$.
The \emph{null line} is the set of all points where the risks of death in the unexposed and the exposed are equal.
It is a diagonal line from the point $(0, 0)$ to the point $(1, 1)$.

The conditional risk of disease given $X = x$ and $C = c$ is
\begin{equation}
  \Pr(D = 1 \given X = x, C = c) = \frac{\Pr(D = 1, X = x, C = c)}{\Pr(X = x, C = c)}.
\end{equation}
The point
\begin{equation}
  p^*_c = \Big(\Pr(D = 1 \given X = 0, C = c), \Pr(D = 1 \given X = 1, C = c)\Big)
\end{equation}
is the \emph{stratum-specific association point} for the subpopulation with $C = c$.
There is one stratum-specific association point for each value of $C$.
Figure~\ref{fig:Rothman_std} shows the stratum-specific association points for the age groups in Table~\ref{tab:strat}.
For both smokers and nonsmokers, the risks of death were much higher in the older age stratum, and these differences in mortality due to age were much larger than the differences due to smoking.

When there is no unmeasured confounding of the causal effect of $X$ on $D$ given $C$, exchangeability and consistency guarantee that
\begin{equation}
  \begin{aligned}
    \Pr(D^x = 1 \given C = c)
    &= \Pr(D^x = 1 \given X = x, C = c) \\
    &= \Pr(D = 1 \given X = x, C = c)
  \end{aligned}
\end{equation}
for $x = 0$ and $x = 1$~\citep{pearl2009causality, hernanrobins2023}.
Thus, the causal and association points are the same (in expectation) for each stratum of $C$.
In our example, the stratum-specific association points represent stratum-specific causal points if stratifying into the 18-64 and 65+ year age groups removes all confounding of the association between smoking and 20-year survival.
In reality, there is likely to be residual confounding by age within these groups.
Better control of confounding could be achieved with a larger number of age groups or a well-chosen regression model.

\section{Standardization and line segments}
If we have two points $p_0 = (x_0, y_0)$ and $p_1 = (x_1, y_1)$ in $\mathbb{R}^2$, we can multiply them by scalars and add them like vectors to get another point in $\mathbb{R}^2$:
\begin{equation}
  a_0 p_0 + a_1 p_1
  = \big(a_0 x_0 + a_1 x_1, a_0 y_0 + a_1 y_1\big).
\end{equation}
As time $t$ goes from zero to one, the point
\begin{equation}
  p(t) = (1 - t) p_0 + t p_1
\end{equation}
goes from $p_0$ to $p_1$ in a straight line.
The line segment connecting $p_0$ and $p_1$ is the set
\begin{equation}
  \big\{(x, y) \in \mathbb{R}^2:
  (x, y) = (1 - t) p_0 + t p_1
  \text{ for some } t \in [0, 1]\big\},
  \label{eq:linesegment}
\end{equation}
which is the set of all points $p(t)$ that we cross while going from $p_0$ to $p_1$ in a straight line.

Because $C$ is not causally affected by $X$, $\Pr(C = 0)$ and $\Pr(C = 1)$ remain the same no matter how we intervene to set $X$.
For any possible value $x$ of $X$,
\begin{equation}
  \begin{aligned}
    \Pr(D^x = 1)
    &= \Pr(D^x = 1 \given C = 0) \Pr(C = 0) \\
    &\qquad + \Pr(D^x = 1 \given C = 1) \Pr(C = 1)
  \end{aligned}
  \label{eq:mcausal}
\end{equation}
by the law of total probability.
If $p_0$ and $p_1$ are the stratum-specific causal points for $C = 0$ and $C = 1$, then
\begin{equation}
  \pmarg = p_0 \Pr(C = 0) + p_1 \Pr(C = 1).
\end{equation}
Because $\Pr(C = 0) + \Pr(C = 1) = 1$, the marginal causal point is on the line segment connecting the stratum-specific causal points.
Other points on this line segment are marginal causal points for populations with the same stratum-specific causal points but different distributions of $C$.
For example, the point $0.5 p_0 + 0.5 p_1$ halfway between $p_0$ and $p_1$ is the marginal causal point for a population where $\Pr(C = 1) = 0.5$.

With association points, line segments correspond to standardized risks of disease~\citep{rothman1975pictorial}.
The standardized risk of disease given $X = x$ is the marginal risk of disease in a hypothetical \emph{standard population} where the prevalence of $C = 1$ equals $\Prst(C = 1)$.
The standardized risk of disease under $X = x$ is
\begin{equation}
  \begin{aligned}
    \Prst(D = 1 \given X = x)
    &= \Pr(D = 1 \given X = x, C = 0) \Prst(C = 0) \\
    &\qquad + \Pr(D = 1 \given X = x, C = 1) \Prst(C = 1).
  \end{aligned}
  \label{eq:std}
\end{equation}
Standardization allows us to calculate marginal risks of disease using the same distribution of $C$ for both exposure groups, which removes confounding by $C$.
The corresponding \emph{standardized association point}
\begin{equation}
  \begin{aligned}
    p^*_\std
    &= \big(\Prst(D = 1 \given X = 0), \Prst(D = 1 \given X = 1)\big) \\
    &= p^*_0 \Prst(C = 0) + p^*_1 \Prst(C = 1)
  \end{aligned}
\end{equation}
is on the line segment connecting the stratum-specific points $p^*_0$ and $p^*_1$.
Each point on this \emph{standardized segment} represents a standardized association point for some standard population.
For example, the point $0.5 p^*_0 + 0.5 p^*_1$ corresponds to a standard population with $\Prst(C = 1) = 0.5$.

If there is no unmeasured confounding of the causal effect of $X$ on $D$ within strata of $C$ and $\Prst(C = 1) = \Pr(C = 1)$, then
\begin{equation}
    p^*_\std
    = p_0 \Pr(C = 0) + p_1 \Pr(C = 1)
    = \pmarg.
\end{equation}
so the standardized point equals the marginal causal point (in expectation).
If we standardize to any other distribution of $C$, we get the marginal causal point for the corresponding standard population.
Table~\ref{tab:standards} shows the three most common standard populations with the causal interpretation of each standardized point in terms of a hypothetical experiment.
In a randomized trial, all three standard populations are the same (up to sampling variation).

\begin{table}
  \centering
  \begin{tabular}{lcl}
    \toprule
    Standard population     & $\Prst(C = 1)$              & Hypothetical experiment \\
    \midrule
    Study sample            & $\Pr(C = 1)$                & \makecell[l]{Study sample under $X = 1$ \\ versus themselves under $X = 0$} \\[10pt]
    Exposed or treated      & $\Pr(C = 1 \given X = 1)$   & \makecell[l]{Exposed versus \\ themselves under $X = 0$} \\[10pt]
    Unexposed or untreated  & $\Pr(C = 1 \given X = 0)$   & \makecell[l]{Unexposed under $X = 1$ \\ versus themselves}\\
    \bottomrule
  \end{tabular}
  \caption{
    Distribution of $C$ for common standard populations and causal interpretation of each in terms of an experiment that compares the risks of disease in same population under two different exposures.
  }
  \label{tab:standards}
\end{table}

In our example, there are 1,072 women aged 18-64 years and 242 women aged 65+ years at the time of the original survey.
Using this age distribution, the standardized risk of death among smokers is
\begin{equation}
  \bigg(\frac{1{,}072}{1{,}314} \times \frac{97}{533}\bigg)
  + \bigg(\frac{242}{1{,}314} \times \frac{42}{49}\bigg)
  \approx 0.306
\end{equation}
and the standardized risk of death among nonsmokers is
\begin{equation}
  \bigg(\frac{1{,}072}{1{,}314} \times \frac{65}{539}\bigg)
  + \bigg(\frac{242}{1{,}314} \times \frac{165}{193}\bigg)
  \approx 0.256.
\end{equation}
This standardized point is plotted in Figure~\ref{fig:Rothman_std} with the label ``All ages''.
If stratifying by age has controlled confounding of the causal effect of smoking on 20-year mortality, the x-coordinate of this point is the counterfactual risk of death in the study sample if we had intervened to make all participants nonsmokers, and its y-coordinate is the counterfactual risk of death if we had intervened (unethically) to make all participants smokers.

Among the 582 smokers in the study sample, there are 533 women aged 18-64 years and 49 women aged 65+ years.
Using this age distribution, the standardized risk of death among smokers is
\begin{equation}
  \bigg(\frac{533}{582} \times \frac{97}{533}\bigg)
  + \bigg(\frac{49}{582} \times \frac{42}{49}\bigg)
  = \frac{97 + 42}{582}
  \approx 0.239,
\end{equation}
which is just the marginal risk of death among smokers in the sample.
The corresponding standardized risk of death among nonsmokers is
\begin{equation}
  \bigg(\frac{533}{582} \times \frac{65}{539}\bigg)
  + \bigg(\frac{49}{582} \times \frac{165}{193}\bigg)
  \approx 0.182.
\end{equation}
This point is plotted in Figure~\ref{fig:Rothman_std}.
It has the same y-coordinate as the crude point because of consistency:
Because we see $D^1$ for smokers, their actual risk of death equals their risk of death under an intervention that makes them smokers.
If stratifying by age has controlled confounding, the x-coordinate of this point is the counterfactual risk of death among smokers had we intervened to make them nonsmokers.
Compared to the standardized point for the study sample, the point for smokers is closer to the stratum-specific point for the 18-64 year age group because smokers were younger on average.

Among the 732 nonsmokers, there are 539 women aged 18-64 years and 193 women aged $\geq 65$ years.
Figure~\ref{fig:Rothman_std} shows the standardized point for this age distribution, which has the same x-coordinate as the crude point because of consistency (i.e., we see $D^0$ for nonsmokers).
If stratifying by age has controlled confounding, then the y-coordinate of this point is the counterfactual risk of death among nonsmokers had we intervened to make them smokers.
Compared to the standardized points for smokers and for the study sample, the point for nonsmokers is closer to the stratum-specific point for the 65+ year age group because nonsmokers were older on average.

\section{Confounding and rectangles}
The set in equation~\eqref{eq:linesegment} is a line segment connecting $p_0$ and $p_1$ because we use the same $t$ for both axes.
If we allow different values of $t$ for the two axes, we get the set
\begin{equation}
  \begin{aligned}
    \big\{(x, y) \in \mathbb{R}^2
    &: x = (1 - t_x) x_0 + t_x x_1
    \text{ and } y = (1 - t_y) y_0 + t_y y_1 \\
    &\quad \text{ for some } t_x, t_y \in [0, 1]\big\},
  \end{aligned}
  \label{eq:rectangle}
\end{equation}
which is a rectangle with sides parallel to the axes and corners at $p_0$ and $p_1$.
It is called the \emph{circumscribing rectangle} because it is the smallest rectangle with sides parallel to the axes that contains these points.
The line segment in equation~\eqref{eq:linesegment} is a diagonal of this rectangle.

The crude association point $\pcrude$ has the same x-coordinate as the standardized point for the unexposed
\begin{equation}
  p_0^* \Pr(C = 0 \given X = 0) + p_1^* \Pr(C = 1 \given X = 0),
\end{equation}
and it has the same y-coordinate as the standardized point for the exposed
\begin{equation}
  p_0^* \Pr(C = 1 \given X = 0) + p_1^* \Pr(C = 1 \given X = 1).
\end{equation}
If the distribution of $C$ is different in the unexposed and the exposed, then $\pcrude$ will be off the standardized segment connecting $p_0^*$ and $p_1^*$ (unless they share an x-coordinate or a y-coordinate).
However, it must be in the circumscribing rectangle with corners at $p_0^*$ and $p_1^*$ and sides parallel to the axes~\citep{rothman1975pictorial}.
This is the \emph{confounding rectangle}, and the standardized segment is one of its diagonals.
The Appendix shows that, in the limit of a large sample, the crude association point $\pcrude$ is on the standardized segment connecting $p_0^*$ and $p_1^*$ if and only if the causal effect of $X$ on $D$ is not confounded by $C$.
Due to sampling variation, this equivalence is only approximate in practice.

\begin{figure}
  \includegraphics[width = \textwidth]{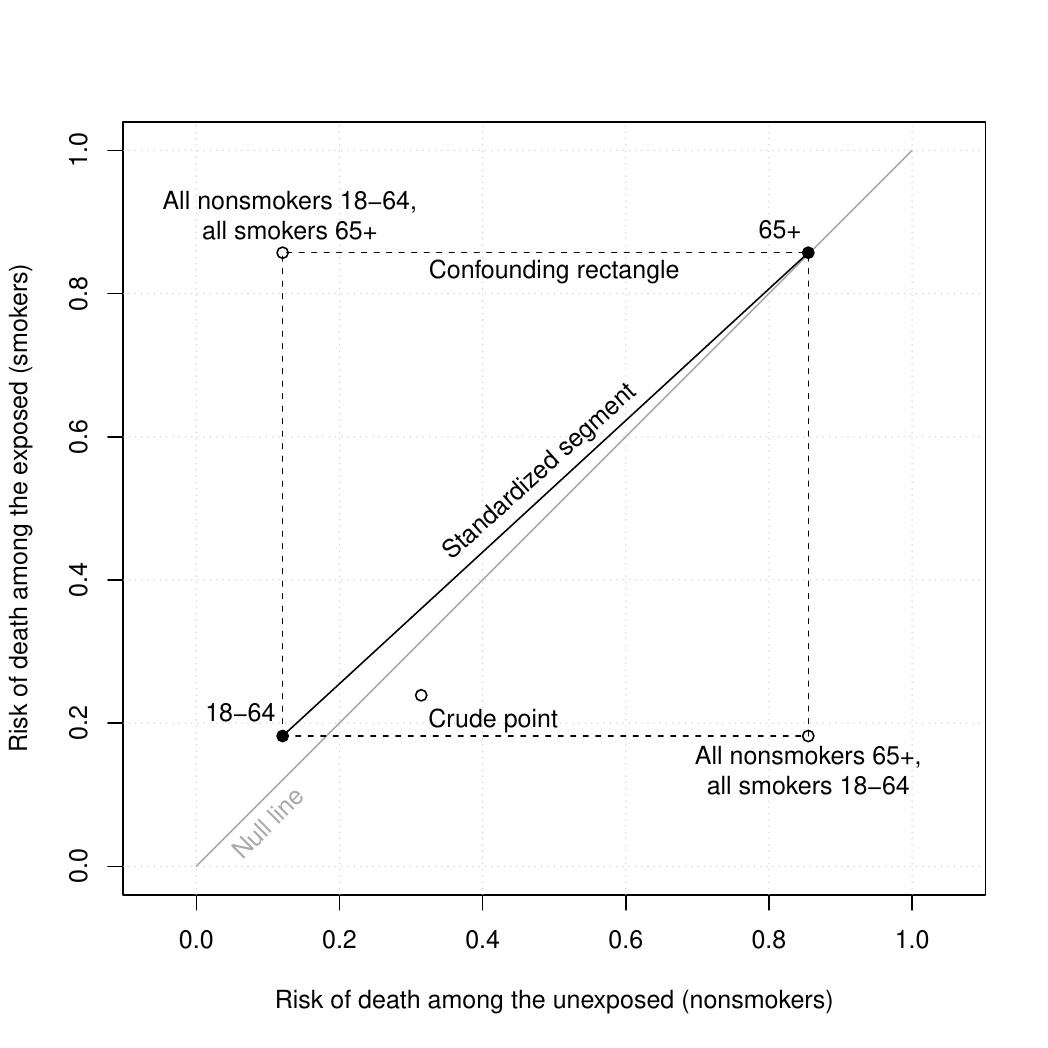}
  \caption{
    A Rothman diagram with the standardized segment (solid) and the confounding rectangle (dashed).
    The stratum-specific points define two corners of the rectangle.
    The other two corners represent the most extreme confounding possible given the stratum-specific risks of death in the exposed and unexposed.
  }
  \label{fig:Rothman_conf}
\end{figure}

Figure~\ref{fig:Rothman_conf} shows a Rothman diagram with the confounding rectangle for the data in Table~\ref{tab:strat}.
The crude point is off the standardized segment because the age distributions for smokers and nonsmokers are different.
The large difference in the risk of death between the two age strata creates a large confounding rectangle that contains points far above and far below the null line.
Thus, unmeasured confounding could create an association that exaggerates or reverses the causal effect of smoking on 20-year mortality.

\section{Confounders with more than two levels}
The generalization of the line segment in equation~\eqref{eq:linesegment} for $k \geq 1$ points $p_1, \ldots, p_k$ is the \emph{convex hull}
\begin{equation}
  \begin{aligned}
    \big\{(x, y) \in \mathbb{R}^2
    &: (x, y) = w_1 p_1 + \ldots + w_k p_k
    \text{ for }  w_1, \ldots w_k \geq 0 \\
    & \quad \text{ such that } w_1 + \ldots + w_k = 1\big\}.
  \end{aligned}
  \label{eq:convexhull}
\end{equation}
The convex hull is the set of all weighted sums of the points $p_1, \ldots, p_k$ where the weights are nonnegative and add up to one.
Geometrically, it is the shape you would get if you stretched a rubber band around all of the points.
It is \emph{convex} because it contains the line segment joining any two of its points.
If $k = 1$, the convex hull is the point $p_1$.
If $k = 2$, the convex hull is the line segment connecting $p_1$ and $p_2$.

Now suppose that $C$ has $k \geq 1$ possible values $c_1, \ldots c_k$.
Then we have $k$ stratum-specific association points $p_1^*, \ldots p_k^*$.
The standardized segment is the convex hull for $k = 2$ strata.
For $k > 2$ strata, we will call the convex hull of the stratum-specific points the \emph{standardized hull}.
Stratum-specific points can be in the interior or on the perimeter of the standardized hull.

\begin{figure}
  \centering
  \includegraphics[width = \textwidth]{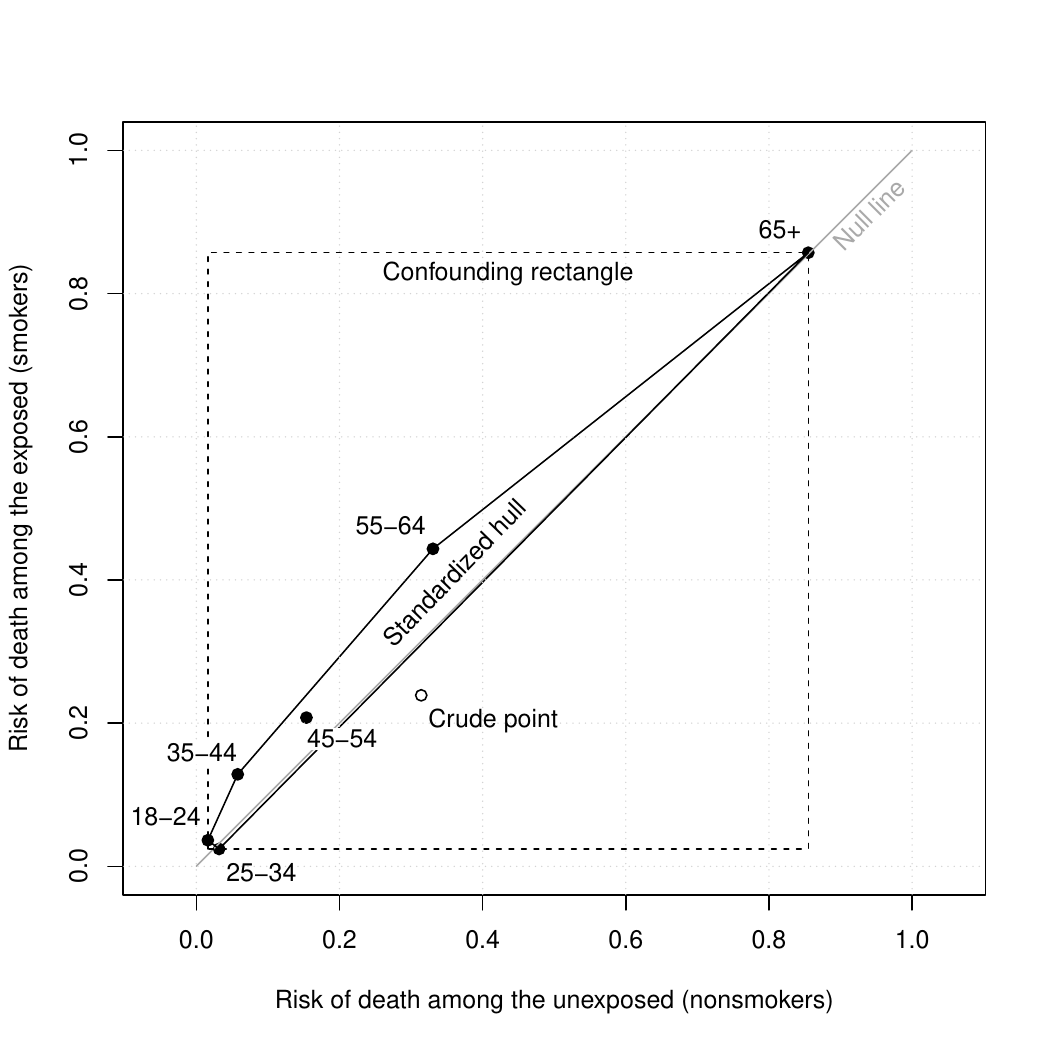}
  \caption{
    The standardized hull (solid) and the confounding rectangle (dashed) for six age strata adapted from~\citet{appleton1996ignoring}.
    Only the top right corner of the confounding rectangle is a stratum-specific point.
  }
  \label{fig:Rothman6}
\end{figure}

Figure~\ref{fig:Rothman6} shows the standardized hull for six age strata adapted from~\citet{appleton1996ignoring}.
The point for the 45-54 year age group is in its interior, but all other stratum-specific points are on its perimeter.
Almost all of the standardized hull is above the null line, where smoking is harmful.

For any $k \geq 1$, we get a confounding rectangle.
Using different sets of weights for the exposed and unexposed to calculate a possible crude point, we can get any point in the circumscribing rectangle for $p_1^*, \ldots, p_k^*$.
It has sides parallel to x-axis with y-coordinates at the smallest and largest stratum-specific risks in the exposed, and it has sides parallel to the y-axis with x-coordinates equal to the smallest and largest stratum-specific risks in the unexposed.
It contains the entire standardized hull.
Unlike the $k = 2$ case, the confounding rectangle for $k > 2$ does not necessarily have stratum-specific points at any of its corners.
Figure~\ref{fig:Rothman6} shows the confounding rectangle for the six age strata.

When there is no confounding of the causal effect of $X$ on $D$ by $C$, the crude point will always be inside the standardized hull.
Each point in the standardized hull is a marginal causal point for some standard population.
For $k > 2$, it is possible to get a crude point inside the standardized hull even when there is confounding.
The Appendix shows that, in the limit of a large sample, there must be confounding if the crude point is outside the standardized hull.
Due to sampling variation, this implication is only approximate in practice.

\section{Effect modification and contour lines}
On a Rothman diagram, every point in the unit square represents a risk in the exposed and a risk in the unexposed.
We can use these risks to calculate measures of association such as the odds ratio, risk ratio, or risk difference.
If we think of the odds ratio as a function
\begin{equation}
  \OR(x, y) = \frac{\sfrac{y}{(1 - y)}}{\sfrac{x}{(1 - x)}},
\end{equation}
then the set of points with $\OR(x, y) = m$ for any given $m$ is called a \emph{contour line} or \emph{contour} of the odds ratio.
Contours can be straight lines or curves, and each possible $m$ corresponds to a different contour line.

We can also find contour lines for the risk ratio $\RR(x, y) = \sfrac{y}{x}$, the risk difference $\RD(x, y) = y - x$, and the hazard ratio.
To estimate the hazard ratio from binary data, we need to assume that there is a constant hazard ratio for the exposed compared to the unexposed during the time interval over which risk is defined.
If so, this hazard ratio is
\begin{equation}
  \HR(x, y) = \frac{\log(1 - y)}{\log(1 - x)}.
\end{equation}
The incidence rate ratio is a special case of the hazard ratio where we assume an exponential distribution of times to disease onset in each exposure group.

\begin{figure}
  \centering
  \includegraphics[width = \textwidth]{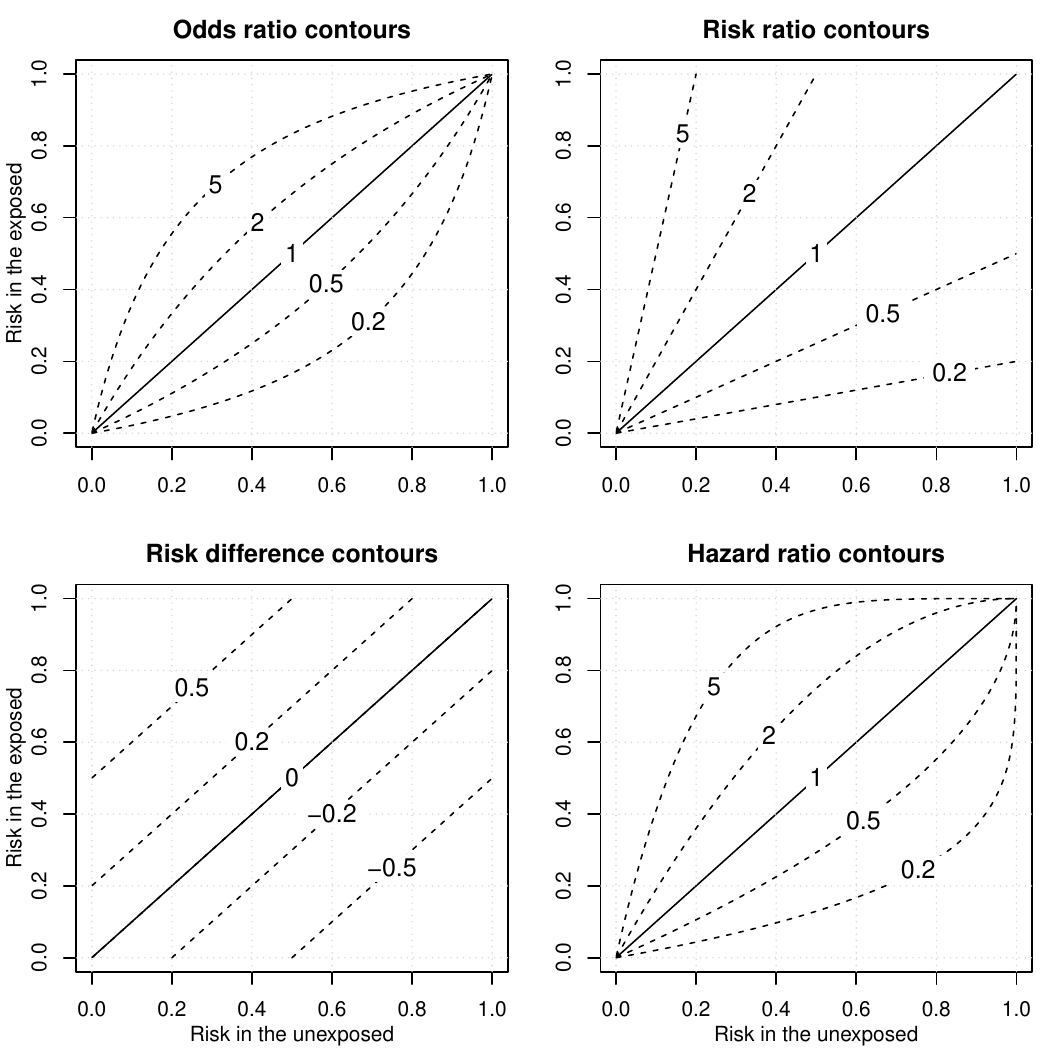}
  \caption{
    Contour lines for the odds ratio, risk ratio, risk difference, and hazard ratio.
    The null line is solid.
    All other contours lines are dashed.
    Each contour line is labeled with the corresponding value of the measure of association.
  }
  \label{fig:contours}
\end{figure}

Figure~\ref{fig:contours} shows contour lines for the odds ratio, risk ratio, risk difference, and hazard ratio~\citep{rothman1975pictorial,richardson2017modeling}.
The null line is a contour line for all measures of association, corresponding to $\RD(x, y) = 0$ and $\OR(x, y) = \RR(x, y) = \HR(x, y) = 1$.
On a Rothman diagram, $C$ is an effect modifier for a measure of association $M$ when the stratum-specific points are on different contour lines of $M$.
If we randomly choose two points in the unit square, it would be an astonishing coincidence if both were on the same contour line of a given measure of association.
Thus, effect modification should be seen as normal, not pathological or unusual.

Table~\ref{tab:effmod} shows stratum-specific measures of association estimated using binomial GLMs that include an interaction between smoking and age group with likelihood-ratio p-values for the null hypothesis that the interaction term coefficient equals zero.
It also shows the common measures of association estimated using the same binomial GLMs without an interaction.
According to the point estimates alone, there is effect modification on all four scales---with a harmful effect of smoking in the 18-64 year age group and a near-null effect of smoking in the 65+ year age group.
However, there is no evidence of effect modification beyond sampling variation on the odds ratio or the risk difference scales.
There is weak evidence of effect modification on the hazard ratio scale and strong evidence on the risk ratio scale.

\begin{table}
  \centering
  \begin{tabular}{lrrrrc}
    \toprule
    Measure of & \multicolumn{3}{c}{Stratum-specific} & Common & 95\% confidence\\
                  \cline{2-4}
    association  & 18-64   & 65+     & p-value  & estimate            & interval \\
    \midrule
    Odds ratio              & 1.622   & 1.018   & 0.353     & 1.537   & (1.119, 2.125) \\
    Risk ratio              & 1.509   & 1.003   & 0.010     & 1.062   & (0.952, 1.166) \\
    Risk difference         & 0.061   & 0.002   & 0.300     & 0.052   & (0.013, 0.091) \\
    Hazard ratio            & 1.563   & 1.008   & 0.085     & 1.316   & (1.034, 1.676) \\
    \bottomrule
  \end{tabular}
  \caption{
    Stratum-specific measures of association for smoking and 20-year mortality with likelihood ratio p-values for the interaction term and common estimates with likelihood ratio 95\% confidence intervals.
  }
  \label{tab:effmod}
\end{table}

\begin{figure}
  \includegraphics[width = \textwidth]{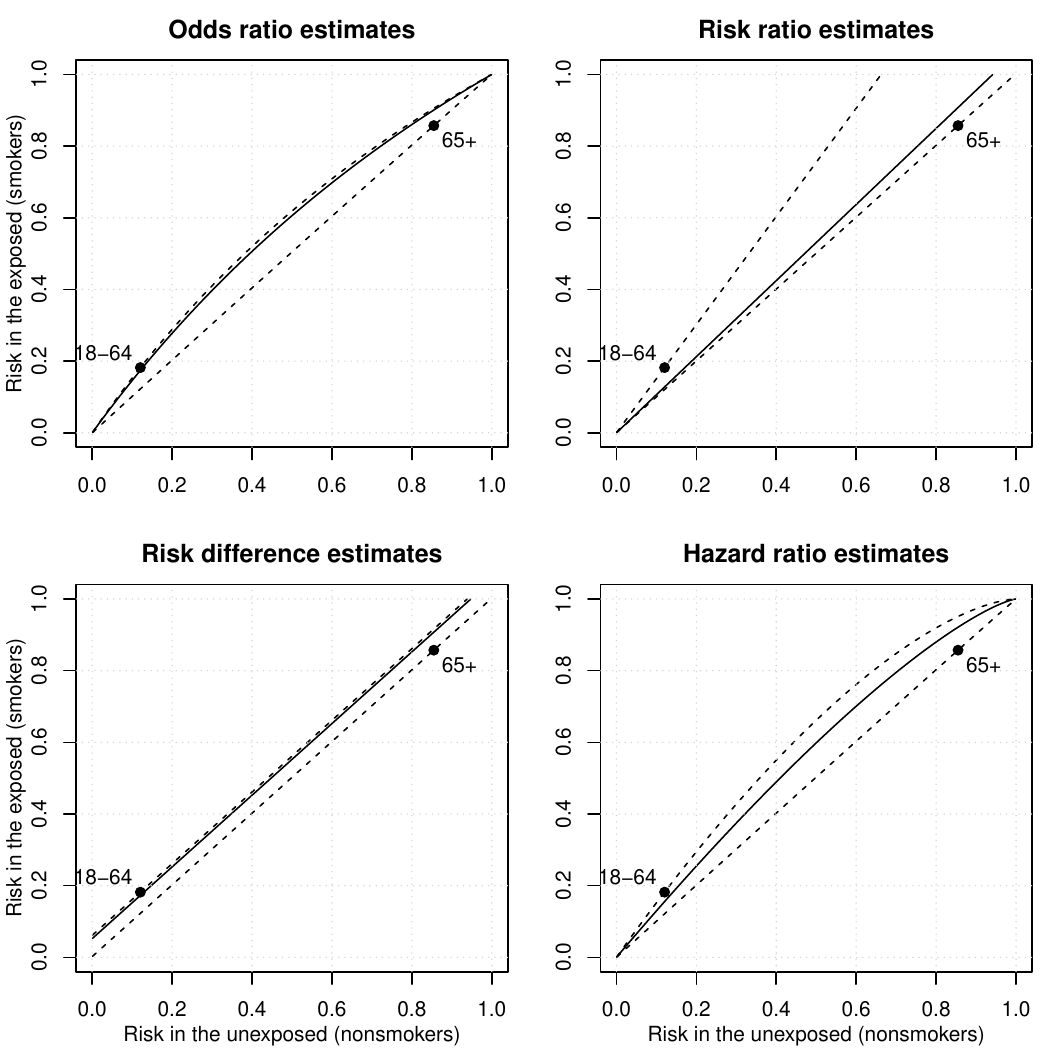}
  \caption{
    Contour lines containing the stratum-specific measures of association (dashed) and common measures of association (solid) on the odds ratio, risk ratio, risk difference, and hazard ratio scales.
  }
  \label{fig:effmod}
\end{figure}

The Rothman diagrams in Figure~\ref{fig:effmod} show the contour lines for the stratum-specific measures of association and for the common measure of association on the odds ratio, risk ratio, risk difference, and hazard ratio scales.
For each measure of association, the points for the 18-64 year age group and the 65+ year age group are on different contour lines, so there is effect modification by age group.
On the odds ratio and risk difference scales, the contour line for the common estimate passes close to both stratum-specific points.
On the risk ratio and hazard ratio scales, the stratum-specific points are farther from the contour line for the common estimate.
These differences partly explain the different p-values in Table~\ref{tab:effmod}.

We can evaluate effect modification whether or not the ``effects'' are causal.
There is effect modification on a given scale when the stratum-specific measures of association are on different contour lines, which can occur whether or not the crude point is on the standardized segment.
There is confounding when the crude point is off the standardized segment, which can occur whether or not the stratum-specific points are on different contour lines for any given measure of association.
Confounding and effect modification are logically independent unless we want to compare stratum-specific causal effects.
To say whether a causal effect differs between strata, we must control confounding first.

\section{Noncollapsibility and curvature}
If $C$ is a risk factor for disease but not a confounder, a measure of association $M$ is said to be \emph{collapsible} if $M = m$ in the marginal table whenever $M = m$ in all strata of $C$~\citep{whittemore1978collapsibility, greenland1999confounding}.
Thus, the marginal and stratum-specific values of a collapsible measure $M$ are equal whenever there is no confounding by $C$ and no effect modification of $M$ by $C$.
The risk ratio and risk difference are collapsible but the odds ratio is not~\citep{samuels1981matching,miettinen1981confounding}.
The hazard ratio is also noncollapsible~\citep{greenland1996absence}.

An equivalent but simpler definition is that $M$ is collapsible if $M = m$ at any standardized point whenever $M = m$ at all stratum-specific points.
Because standardization removes confounding by $C$, this is equivalent to the definition that assumes no confounding by $C$.
Leaving confounding out of the definition of collapsibility avoids a common source of confusion in epidemiology and biostatistics, where students are often taught to detect confounding by looking for a change in an estimated measure of association upon adjustment for a covariate.
This ``change-in-estimates'' method can spuriously detect confounding when it is based on a noncollapsible measure of association~\citep{miettinen1981confounding, greenland1999confounding}.

\begin{figure}
  \centering
  \includegraphics[width = \textwidth]{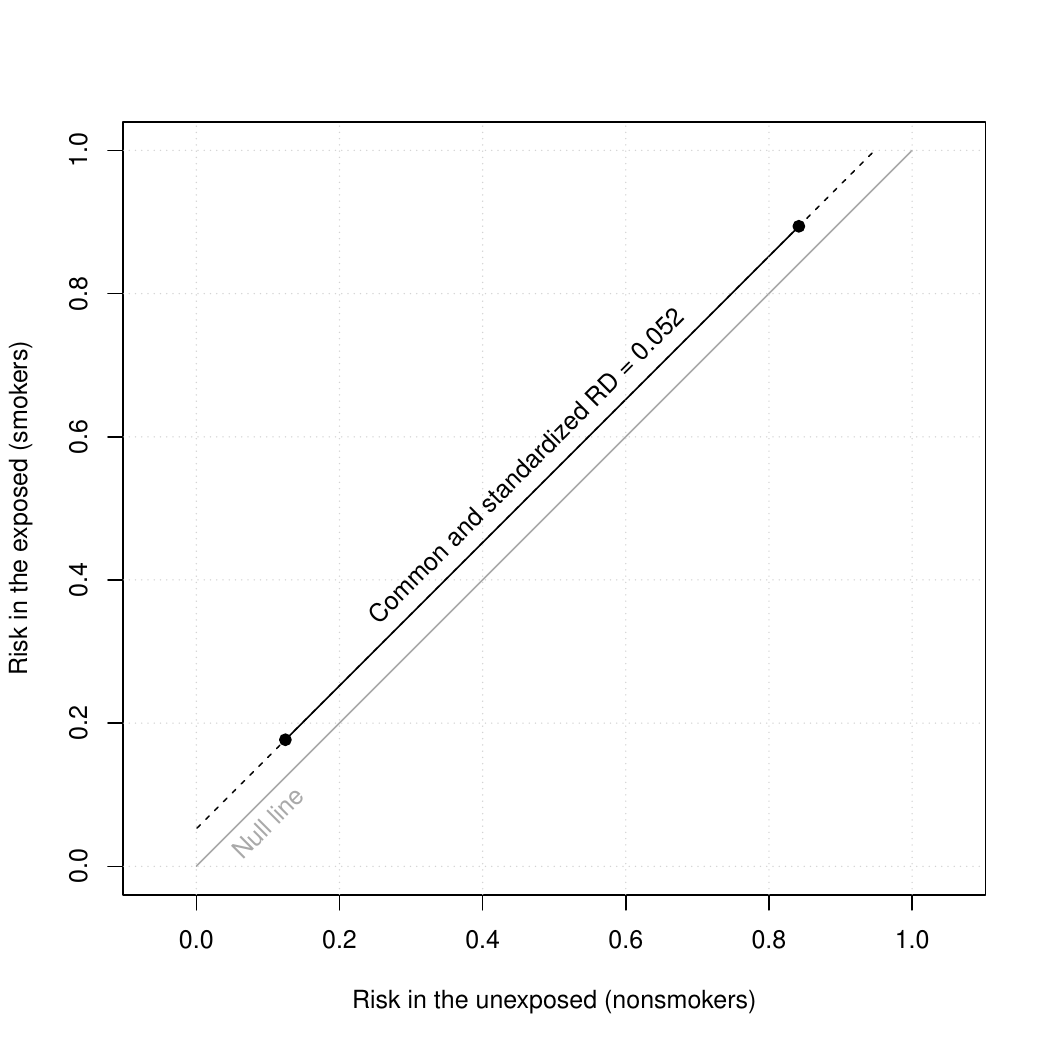}
  \caption{
    Collapsibility along a straight contour line.
    The entire standardized segment is contained on the common risk difference contour.
  }
  \label{fig:collapsible}
\end{figure}

\begin{figure}
  \centering
  \includegraphics[width = \textwidth]{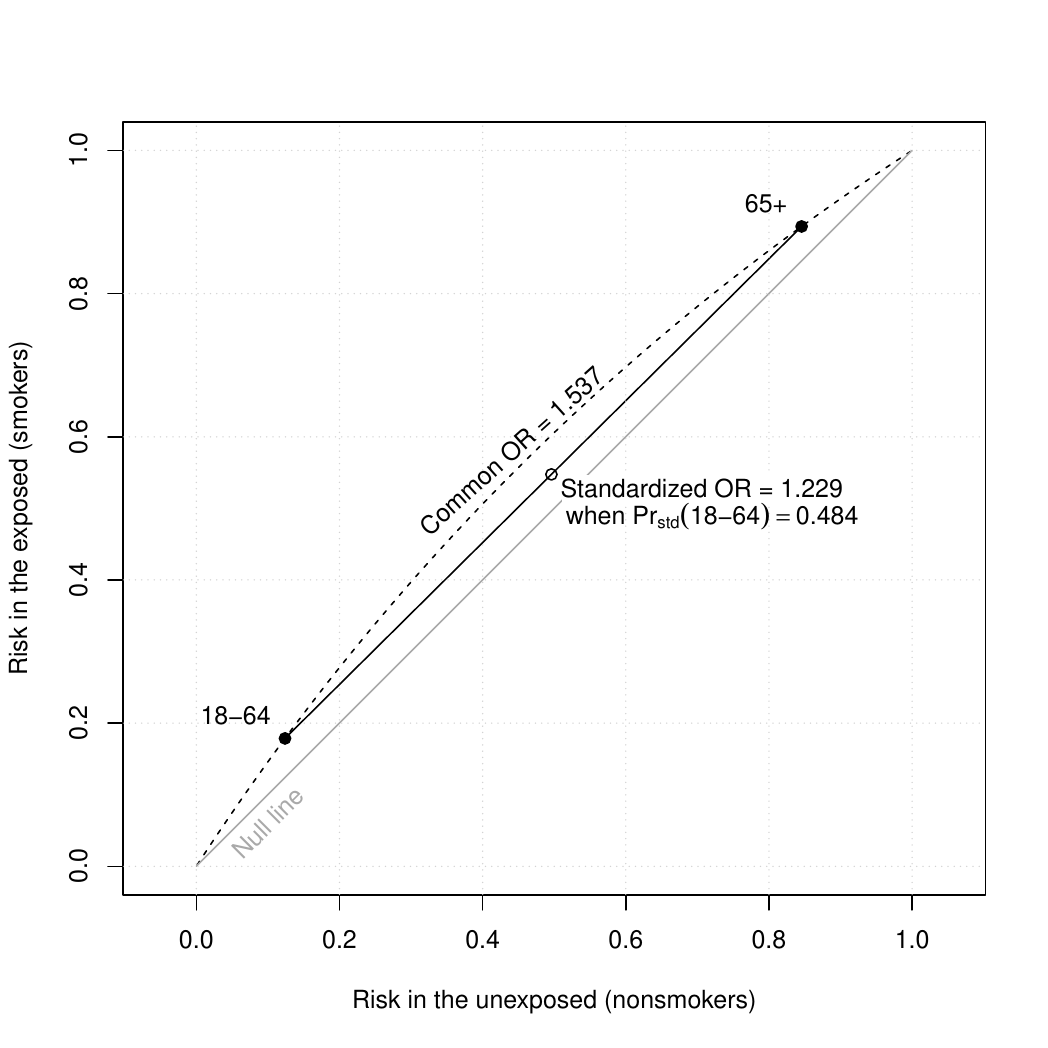}
  \caption{
    Noncollapsibility along a curved contour line.
    Except for its endpoints at the stratum-specific points, the entire standardized segment is between the null line and the contour for the common odds ratio.
  }
  \label{fig:noncollapsible}
\end{figure}

A measure $M$ is collapsible if and only if all of its contour lines are straight.
In Figure~\ref{fig:contours}, the collapsible risk ratio and risk difference have straight contours but the noncollapsible odds ratio and hazard ratio have curved contours.
If all contours of $M$ are straight and all stratum-specific points are on a contour line where $M = m$, then the standardized hull collapses to a line segment along this contour, so $M = m$ at all standardized points.
Since this is true for any $m$, the measure $M$ is collapsible.
This situation is illustrated in Figure~\ref{fig:collapsible}, where all standardized risk differences equal the stratum-specific risk difference of $0.052$.
If any contour $M = m$ is curved, we can find points $p_0$ and $p_1$ on this contour whose standardized segment contains points where $M \neq m$.
Thus, $M$ is not collapsible.
This situation is illustrated in Figure~\ref{fig:noncollapsible}, where the stratum-specific odds ratio is $1.537$ but the standardized segment has a minimum odds ratio of $1.229$, which occurs in a standard population that is approximately $48.4\%$ aged 18-64 years and $51.6\%$ aged $\geq 65$ years.
Because the null line is a contour for all measures of association and it is straight, all measures of association are collapsible under the null hypothesis.

A change in an estimated measure of association upon covariate adjustment is caused by a combination of confounding, effect modification, and noncollapsibility.
Generating a large difference between estimated measures through noncollapsibility alone requires a common stratum-specific estimate far from the null (so the contour has high curvature) and large differences in risk between strata.
Also, noncollapsibility generates bias toward the null but never across the null.
Confounding and effect modification generate large biases---including bias across the null---much more easily, and they affect all measures of association.
Effect modification is a much more serious threat than noncollapsibility to the change-in-estimates method of detecting confounding.

\section{Geometry before regression}
The geometric perspective pioneered by~\citet{rothman1975pictorial} and elaborated above gives standardization its proper place as the bedrock of causal inference and makes clear distinctions between confounding, effect modification, and noncollapsibility.
In the minds of many epidemiologists, these concepts are almost inextricably entangled.
One reason for this is that they are often taught to look for confounding and causal effects in the coefficients of regression models.
The change-in-estimates method for detecting confounding encourages conflation of noncollapsibility and confounding while concealing the danger of ignoring effect modification.
Effect modification and confounding are bound together by the fact that a regression model may need to include interaction terms in order to control confounding in the presence of effect modification.
Worst of all, regression has displaced standardization in epidemiologic methods training even though they are more useful together than either is by itself.

Regression models are wonderful statistical tools, but they were never intended to teach basic concepts in causal inference.
When used for this purpose, they generate confusion and mythology.
They should be introduced only after students have achieved a clear understanding of standardization, confounding, effect modification, and collapsibility.
Because people usually think more clearly in terms of pictures than mathematical symbols, Rothman diagrams are an excellent way to introduce these concepts in an intuitive and unified framework.
As it was written over the entrance to Plato's Academy, so it can be said of regression in causal inference: ``Let no one ignorant of geometry enter here.''

\section*{Supplementary material}
The file ``GCIepi.R'' contains code in \texttt{R} to produce Tables~1-3 and 5 and Figures 1-7 as well as the 2x2 tables for all six age groups that appear in Figure~\ref{fig:Rothman6}.
Instructions for running it are given in comments near the beginning of the file.

\section*{Acknowledgements}
I want to to thank the students of STA 6177/PHC 6937 (Applied Survival Analysis) at the University of Florida in 2014-2016 and PUBHEPI 8430 (Epidemiology 4) at Ohio State in 2019-2023 for working with early versions of this material.
I also want to thank M. Elizabeth Halloran, Wasiur R. KhudaBukhsh, Bo Lu, Nick Mandarano, Kunjal Patel, Micaela Richter, Kenneth Rothman, and Patrick Schnell for their useful comments.
This work was supported by National Institute of Allergy and Infectious Diseases (NIAID) grant R01 AI116770 (PI: Eben Kenah) and National Institute of General Medical Sciences (NIGMS) grant U54 GM111274 (PI: M. Elizabeth Halloran).
The content is solely the responsibility of the author and does not represent the official views or policies of NIAID, NIGMS, or the National Institutes of Health.

\bibliographystyle{chicago}
\bibliography{GCIepidemiology}

\clearpage
\appendix
\section*{Appendix: Confounding and crude points}
When the causal effect of $X$ on $D$ is confounded by $C$, the crude association point $\pcrude$ must be off the standardized segment connecting $p_0^*$ and $p_1^*$ (unless they share an x-coordinate or a y-coordinate)
In a causal DAG that contains all three variables, there is an open backdoor path from $X$ to $D$ that goes through $C$~\citep{greenland1999causal}.
Because of the open path from $C$ to $D$ that does not go through $X$, at least one of the conditional risks of disease given $X$ differs between strata of $C$, so the stratum-specific association points $p^*_0$ and $p^*_1$ are different.
Because of the open path from $C$ to $X$, the conditional distribution of $C$ given $X$ is different when $X = 0$ and $X = 1$.
Therefore, $p_0^*$ and $p_1^*$ define a line segment that does not contain $\pcrude$.
If $p_0^*$ and $p_1^*$ have the same x-coordinate (i.e., risk among the unexposed is independent of $C$) or the same y-coordinate (i.e., risk among the exposed is independent of $C$), then the confounding rectangle and the standardized segment are the same.

When there is no confounding of the causal effect of $X$ and $D$ by $C$, the crude point must be on the standardized segment.
If there is no open path from $C$ to $D$, then the stratum-specific association points $p^*_0$ and $p^*_1$ are the same, so the confounding rectangle and the standardized segment both collapse to that point.
If there is no open path from $C$ to $X$, then the distribution of $C$ is the same for $X = 1$ and $X = 0$, so $\pcrude$ is on the standardized segment.
Therefore, we have either $\pcrude = p_0^* = p_1^*$ or $\pcrude$ is on the standardized segment connecting $p^*_0$ and $p^*_1$.

Due to sampling variation, we may get a $\pcrude$ that is on the standardized segment despite confounding or a $\pcrude$ that is off the standardized segment despite no confounding.
This sampling variation disappears as the sample size $n \rightarrow \infty$, so the equivalence between confounding and the crude point being off standardized segment holds exactly in the limit of a large sample from an infinite population.
In finite samples, this equivalence holds only approximately.

For $k > 2$, sampling variation can produce a crude point outside the standardized hull when there is no confounding.
This disappears in the large-sample limit as $n \rightarrow \infty$, so a crude point outside the standardized hull implies confounding.
In finite samples, this implication holds only approximately.

\end{document}